\documentclass[12pt]{article}
\usepackage{amsmath}
\usepackage{amssymb}
\usepackage{hhline}
\usepackage{array}
\usepackage[dvips]{graphicx}
\usepackage{graphicx}
\usepackage{epsfig}
\usepackage{color}
\usepackage{bbm}
\newtheorem{theorem}{\bf Theorem}[section]
\newtheorem{remark}[theorem]{\bf Remark}
\newtheorem{assumption}[theorem]{\bf Assumption}

\newtheorem{corollary}[theorem]{\bf Corollary}

\newtheorem{proposition}[theorem]{\bf Proposition}

\newenvironment{Proofc}[1]{\smallskip\par\noindent\textsc{#1}\quad}%
  {\hfill$\Box$\bigskip\par}

\numberwithin{equation}{section}

\title{A priori gradient bounds for fully nonlinear parabolic equations and applications to porous medium models}
\author{H. Hajj Chehade$^{1,2}$, M. Jazar$^{2}$ and R. Monneau$^3$}

\begin{document}

\maketitle

\footnotetext[1]{\noindent LAMFA, University of Picardie Jules Vernes, UFR des Sciences
33, rue Saint-Leu 80039 Amiens Cedex 1, France.
E-mail: hana.hajj.chehade@u-picardie.fr (H. Hajj Chehade)}
\footnotetext[2]{\noindent LaMA-Liban, Azm Research Center, Lebanese University, Tripoli, Lebanon.
E-mail: mjazar@ul.edu.lb (M. Jazar)}
\footnotetext[3]{\noindent CERMICS, Ecole des Ponts ParisTech, 6 et 8  avenue Blaise Pascal, Cit\'{e} Descartes, 77455 Marne-La-Vall\'{e}e Cedex 2, France.
E-mail: monneau@cermics.enpc.fr (R. Monneau)}
\noindent\footnotesize .\\

\noindent\small{\textbf{Abstract}. We prove  a priori  gradient bounds for
classical solutions of the fully nonlinear
 parabolic equation
  $$u_{t}=F(D^2u ,D u,u,x,t).$$ Several applications  are given,  including
 the standard porous medium equation.
}
\noindent{\small\textbf{AMS Subject Classifications:}}
{\small{35B50, 35D40, 35K55,  35K65, 36S05, 74G45, 86A05}}\hfill
\medskip


\section{Introduction}
\noindent  Consider the general fully nonlinear parabolic problem
\begin{equation}\label{pborigine}
u_{t}=F(D^2u ,D u,u,x,t),\quad (x,t)\in Q:={\mathbb{T}}^d\times(0,+\infty),
\end{equation}
\begin{equation}\label{initdata}
u(x,0)=u_{0}(x),\quad x\in {\mathbb{T}}^d.
\end{equation}
To simplify our arguments, we consider the case of the $d$-dimensional torus ${\mathbb{T}}^d:=(\mathbb{R}/{\mathbb{Z}})^d$ for $d\ge1.$ Note that up to the price of technicalities, the case of the whole space ${\mathbb{R}}^d$ could be also considered.
The aim of the paper is to find assumptions on $F$ in order to obtain, for all $t\ge 0$, a gradient bound on $D u $ of the form
\begin{equation}\label{boundongradient}\|D u(\cdot,t)\|_{{\infty}}\le\|Du_0\|_{{\infty}}.\end{equation}
As an application of our general approach, we prove  gradient estimate (\ref{boundongradient})   for  the weak nonnegative solution of the standard porous medium equation
\begin{equation}\label{PME}
u_{t}=\Delta u^m,\quad (x,t)\in Q,
\end{equation}   where  $1\le m\le 1+\frac{4}{3+d}.$ For such range of $m$, this result is new.

\noindent Similar gradient estimates are  given for the  problem

 \begin{equation}\label{PMEGprob}
u_{t}=\Delta G(u),\quad(x,t)\in Q,
\end{equation}
for some class of functions $G$, and for the doubly nonlinear problem
 \begin{equation}\label{psimodelegeneral}
u_{t}=\textrm{div}\,\left(\psi(u,|Du|^2)Du\right),\quad(x,t)\in Q,
\end{equation}
for some class of functions $\psi$. Our estimate will be applied to two specific examples of equation (\ref{psimodelegeneral})  arising in hydrology (and this was our initial motivation for this work). These two examples are the following equations

\begin{equation}\label{modelehydrologie}
u_{t}=\textrm{div}\,\left(u(1-u)\frac{Du}{1+|Du|^{2}}\right),\quad (x,t)\in Q,
\end{equation}

\noindent and
\begin{equation}\label{modelehydrologiesimplifie}
u_{t}=\textrm{div}\,\left(u(1-u)Du\right),\quad (x,t)\in Q.
\end{equation}
Equation (\ref{modelehydrologiesimplifie}) derives from equation (\ref{modelehydrologie}) as an approximation for small gradients. In these two equations the function $u$ represents the height of
 the sharp interface between salt and fresh water in a confined aquifer, see for instance \cite{Bear, DH, Josselin De Jong}.

\subsection{Main results}
\noindent In this subsection, we will present our main results. To this end, we will start by an assumption on the function $F$ appearing in equation (\ref{pborigine}). In order to write this assumption, we need to introduce some notation.

 For two symmetric matrices  $X=(x_{ij})_{1\leq i,j\leq d}$ and $Y=(y_{ij})_{1\leq i,j\leq d}$ in $\mathbb{R}^{d\times d}$,
we denote by  $X:Y$ the inner scalar product $\underset{i,j=1,..,d}{\sum}x_{ij}y_{ji}=tr(XY).$ Moreover for $p,q\in{\mathbb{R}}^d$, we set $(X\cdot p)_i=\underset{j=1,..,d}{\sum}x_{ij}p_{j}$ and $p\cdot q=\underset{j=1,..,d}{\sum}p_{j}q_{j}$. For later use, we also denote by $tr (X)$ the trace of $X$.

\begin{assumption}\label{assumpF}\textbf{ }\\
Let $d\ge 1$ and let $\underline u<\overline u$ and $M>0$ be three real numbers,
 $$E:= S^d\times\overline{B(0,M)}\times[\underline{u},\overline{u}]\times Q\quad
\subset \quad S^d\times {\mathbb{R}}^{d}\times {\mathbb{R}}\times Q,$$  where $S^d$ is the
  space of $d\times d$-symmetric real matrices and $F=F(X,p,u,x,t)$ is a real function  defined on $E$
satisfying the following conditions:
\noindent\item[i)]\textbf{\textit{Regularity}}: $F\in C^{1}(E).$
\noindent\item[ii)]\textbf{\textit{Degenerate parabolicity}}:
   For all $X,Y\in S^d$, $(p,u,x,t)\in {\mathbb{R}}^d\times \mathbb{R}\times Q$, we have the implication
  $$\mbox{if } X\le Y, \quad \mbox{ then } F(X,p,u,x,t)\le F(Y,p,u,x,t).$$
\noindent\item[iii)]\textbf{\textit{Differential inequality}}:\begin{equation}\label{A}
-D_XF:X^2+|p|^2D_uF+D_x F\cdot p\le 0,
\end{equation}
for all $(X,p,u,x,t)\in E$ such that $ X\cdot p=0$.
\end{assumption}

 \begin{theorem}(\textbf{{A priori gradient bound for  Problem (\ref{pborigine})}})
 \label{general case}\\
Suppose that $F$ satisfies Assumption \ref{assumpF} and let us assume the existence of a function  $u\in C^{3}(\overline{Q})$ solution of the
 parabolic problem (\ref{pborigine}), (\ref{initdata}) such that for all $t\ge0$, we have
\begin{equation}\label{$H_{u}$}
\underline{u}\le u(\cdot,t)\le \overline{u} \quad \quad\mbox{ and } \quad\quad\|Du(\cdot,t)\|_{L^{\infty}({\mathbb{T}}^d)}\le M\mbox{;}
 \end{equation}
then for all $t\ge 0$
\begin{equation}\label{uniformboundonnabla}
\|Du(\cdot,t)\|_{L^{\infty}({\mathbb{T}}^d)}\le \|Du_{0}\|_{L^{\infty}({\mathbb{T}}^d)}.
\end{equation}
Moreover, if $F$ satisfies the following condition
\begin{equation}\label{constantsolutions}
F(0,0,C,x,t)=0, \mbox{ for any constant } C\in[\underline u,\overline u];
\end{equation} then for all $t\ge 0$
\begin{equation}\label{uniformboundonu}
\underset{{\mathbb{T}}^d}{\min}\, u_0\le u(\cdot,t) \le\underset{{\mathbb{T}}^d}\max \,u_0.
\end{equation}
\end{theorem}

\begin{remark}(\textbf{Generalization})\\
Up to adapt the proofs with certain technicalities, it would also be possible to get a similar result for the same equation on $\Omega \times (0,T)$ with $\Omega={\mathbb{R}}^d $ or ${\mathbb{T}}^d$, and $T>0$.
\end{remark}

\begin{corollary}\label{PMEGresult}(\textbf{A priori gradient  bounds for Problem (\ref{PMEGprob})})\\
Let $d\ge 1$ and $G\in C^{3}([\underline{u}, \overline{u}])$, satisfying
$G'\ge 0 \mbox{ on } [\underline{u}, \overline{u}]$  and
\begin{equation}\label{Gcondition}
\frac{d-1}{4}(G'')^{2}\le-G^{(3)}G'\quad\mbox{ on }[\underline{u}, \overline{u}].\end{equation}
\noindent Let us assume the existence of $u\in C^{3}(\overline{Q})$ solution of
(\ref{PMEGprob}), (\ref{initdata}) such that $u$ and $Du$  satisfy condition (\ref{$H_{u}$}) for some $M>0$;
then for all $t\ge0$, $u$ and $Du$ satisfy a priori bounds (\ref{uniformboundonu}) and (\ref{uniformboundonnabla}).
\end{corollary}
The next corollary gives a new result for the standard porous medium equation (\ref{PME}).

\begin{corollary}\label{PMEresult}(\textbf{Application to the standard porous medium  equation})\\
Let $d\ge 1$,  $\displaystyle{1\le m\le1+\frac{4}{3+d}}$ and $u_0\in W^{1,\infty}({\mathbb{T}}^d)$ with $u_0$ nonnegative. Then there exists a
 unique nonnegative weak solution $u\in L^{\infty}(0,+\infty; W^{1,\infty}({\mathbb{T}}^d)) \cap C([0,\infty);L^{1}({\mathbb{T}}^d))$ of
Problem (\ref{PME}), (\ref{initdata}). Moreover, $u$
satisfies   bounds (\ref{uniformboundonu}) and (\ref{uniformboundonnabla}).
\end{corollary}

 \subsection{Brief review of the literature}
 When $F$ does not depend  on $u$, certain gradient estimates of the form (\ref{boundongradient}) do exist in the literature  for solutions of Problem (\ref{pborigine}), (\ref{initdata}). For instance (see \cite{BL}), such estimates hold true for solutions of the equation
$$u_t=\textrm{div} (\psi_0(|D u|^2)D u),$$
 in any dimension $d\ge 1$ and under some assumptions on $\psi_0.$ Notice also that Assumption \ref{assumpF} iii) is always satisfied if $F$ does not depend on $u$ neither on $x.$

\noindent Let us notice that in \cite{lions}, a uniform gradient estimate is given for bounded solutions of the elliptic equation
$$\lambda u +F_0(D^2u,Du,x)=0,\quad x\in {\mathbb{R}}^{d},$$
where $\lambda>0$ and under some hypothesis on $F_0$. Such elliptic equation can be seen as a time discretization of evolution equation of type (\ref{pborigine}) when $F$ does not depend on $u.$\\
In \cite{MaPa},  a uniform gradient bound is given  for  the general quasilinear equation $$u_t=a(x,t,u,Du)+\underset{i,j=1}{\overset{d}{\sum}}a_{ij}(x,t,u,Du)\partial_{x_ix_j}u,\quad(x,t)\in \mathcal{O}\times(0,T) , $$ where  $\mathcal{O}$ is an open bounded domain of ${\mathbb{R}}^{d}$ and  $T$ is small enough. This gradient bound depends both on the gradient of the initial data and on a bound on the gradient of the boundary of $\mathcal{O}.$

In \cite{Lu380}, a bound on $D (G^{\alpha}( u))$ for some $\alpha>\frac12$ is given for   solutions of (\ref{PMEGprob}), (\ref{initdata}). The bound is given  by a large constant depending on the initial data.\\
Applying this result to the standard porous medium equation (\ref{PME}), (\ref{initdata}), with the choice  $G(u)=u^{m}$ for $m$ in a certain range, this gives for all $t\ge 0$
\begin{equation}\label{boundongradienttype2}
\|D u(\cdot,t)\|_{\infty}\le M_0(\|Du_0\|_{\infty}),\end{equation}
where
 $M_0(\|Du_0\|_{\infty})$ is a large constant depending on the initial data.\\
 Notice that a  particular application of the general results in  \cite{JL} also gives a similar bound for different range of values of $m$.\\
For the same equation, in any dimension $d\ge 1$, it is possible to deduce from
\cite{MP392} certain interior gradient estimates assuming some local integrability of the gradient of a power of the solution.

\noindent For the porous medium  equation, the only known gradient bound of  the form (\ref{boundongradient}) is given in   \cite[Proposition 15.4, p. 359]{vasquez}  where the author proves the general result  in dimension $d=1$ and for $t\ge 0$
$$\|Du^{m-1}(\cdot,t)\|_{\infty}\leq \|D u^{m-1}_0\|_{\infty}.$$
In particular, for $m=2$, this gives an estimate of the form (\ref{boundongradient}). Notice that our corollary (\ref{PMEresult}) gives gradient estimate of the form (\ref{boundongradient}) at least in dimension $d=1$ for $1\le m\le 2.$

We underline the fact that in this paper, we focus on gradient estimates of the form (\ref{boundongradient}) and not of the less precise form (\ref{boundongradienttype2}).



\subsection{Organization of the paper}
\noindent In section 2, we prove the main results. Section 3 is dedicated
 to  equation (\ref{psimodelegeneral}) and hydrological models (\ref{modelehydrologie}) and (\ref{modelehydrologiesimplifie}).

\section{Proof of the main results}
\noindent Before proving Theorem \ref{general case}, we start by some notation
and recall the  maximum principle for solutions $u$ of Problem (\ref{pborigine}).

\bigskip
\noindent\textbf{Notation}\\
\noindent\textit{i) For any  $i,j,k=1,..,d$, we  consider the following derivatives $$\begin{array}{lll}
 \displaystyle{u_{i}=D_{x_i}u=\frac{\partial u}{\partial x_i}}&\quad
 \displaystyle{ u_{it}=D_tu_{i}=\frac{\partial^2 u}{\partial t\partial x_i}},&
\quad \displaystyle{ u_{ij}=D^{2}_{x_ix_j}u=\frac{\partial^2 u}{\partial x_i\partial x_j.}}
\end{array}$$
\noindent ii) For $T>0$, we denote by $C^{2,1}({\mathbb{T}}^d\times [0,T])$ the space of continuous functions such that the derivatives
 $u_{t}$, $u_{i}$  and $u_{ij}$    exist and are continuous on ${\mathbb{T}}^d\times [0,T]$ for $i,j=1,\ldots, d$.
}
\begin{proposition}\label{comparisononu}(\textbf{{Maximum principle for solutions of $(\ref{pborigine})$}})\\
Let $T>0$ and $F$ be a function satisfying conditions i) and ii) of Assumption \ref{assumpF}.
   Assume the existence of $u$ and $v \in C^{2,1}({\mathbb{T}}^d\times (0,T])\cap C({\mathbb{T}}^d\times [0,T])$ two solutions
 of equation (\ref{pborigine})  satisfying condition
  (\ref{$H_{u}$}).    If  $u(\cdot,0)\leq v(\cdot,0)$, then $u(\cdot,t)\leq v(\cdot,t)$ for all $t\in[0,T]$.
\end{proposition}

\noindent This result  seems standard, even if we have no references for this specific result. For sake of convenience,
a short proof is given in the appendix.

\begin{Proofc}{\textbf{Proof of  Theorem \ref{general case}.}}

\noindent Bounds (\ref{uniformboundonu}) on $u$ are  a direct application of Proposition \ref{comparisononu} and condition (\ref{constantsolutions}), comparing the solutions with both the minimum and the maximum of $u_0$ which are two other constant solutions  of the same equation.\\
The uniform bound on $Du$ will be done in the following four steps:
\bigskip

\noindent\textbf{Step 1: Differential equation for $|Du|^{2}$ }\\
\noindent Differentiating  equation (\ref{pborigine}) yields
\begin{equation}
u_{it}=D_XF: D^2u_i+ D_pF\cdot Du_i+u_iD_uF+D_{x_i}F.
\end{equation}
Multiplying by $u_i$, we get
\begin{equation}\label{foisu_i}
 \displaystyle{\left(\frac{u_i^{2}}{2}\right)_t=D_XF:u_i D^2u_i+D_pF\cdot u_i Du_i+u_i^2D_uF+u_iD_{x_i}F}.
\end{equation}
Since
\begin{equation}\label{Egalite 1}
D(|Du|^{2}/2)=D^{2}u\cdot Du
\end{equation}
and
\begin{equation}\label{Egalite 2}
D^{2}({|Du|^{2}/2})=\sum_{i} u_i D^{2} u_i+(D^{2}u)^{2},
\end{equation}
setting
 $w:=|Du|^{2}/2$, by summation on $i$ in (\ref{foisu_i}), we obtain
\begin{eqnarray}\label{equationenw}
w_t=D_XF\cdot( D^2w-(D^2u)^2)+D_pF\cdot D w+2wD_uF+D_xF\cdot D u,
\end{eqnarray}
where we used equality (\ref{Egalite 2}) for the first  term and (\ref{Egalite 1}) for the second one.
\bigskip

\noindent\textbf{Step 2: Differential equation at the maximum of  $|Du|^2$}\\
Now, for all $t\ge 0$, set $\displaystyle M(t):=\max_{{\mathbb{T}}^d} w(\cdot,t)$.
Let  $t_{0}>0$ and $(x_{0},t_{0})$  be a point such that $w(x_{0},t_{0})=M(t_{0})$. At $(x_0,t_0)$, we have
\begin{equation}\label{x_0,t_0,D}Dw(x_0,t_0)=0,
\end{equation}
and
\begin{equation}\label{x_0,t_0,D^2} D^2w(x_0,t_0)\le 0.
\end{equation}
From Assumption \ref{assumpF} ii), we deduce that $D_{X}F\ge 0$
and then
 $D_XF:(D^2w(x_0,t_0))\le 0$ by (\ref{x_0,t_0,D^2}).
Similarly, we get
$D_pF\cdot Dw(x_0,t_0)=0$ by (\ref{x_0,t_0,D}). This implies that
\begin{equation}\label{w|_max}
w_t\leq -D_XF:(D^2u)^2+2wD_uF+D_x F\cdot D u.\end{equation}

\noindent\textbf{Step 3:  Inequality  $M'\le 0$}\\
 In this step, we prove that, in the viscosity sense, for all $t>0$, we have
\begin{equation}\label{viscosity}M'(t)\le0.\end{equation}

\noindent Let $V$ be a neighborhood of $t_0\in(0,+\infty)$ and  $\phi\in C^{1}(0,+\infty)$
 verifying
$$\left\{\begin{array}{ll}
M(t)\le \phi(t) & \mbox{on V,}\\
M(t_{0})=\phi(t_0);&
\end{array}\right.$$
then $w(x_{0},t)\le\phi(t)$
and $w(x_{0},t_{0})=\phi(t_{0})$,
which implies that $\phi'(t_{0})=w_{t}(x_{0},t_{0})$.
 Note that (\ref{x_0,t_0,D}) and (\ref{Egalite 1}) imply that $X\cdot p=0$ with $p=Du(x_{0},t_{0})$ and $X=D^2u(x_{0},t_{0})$.
 Therefore, we deduce from (\ref{w|_max}) and Assumption \ref{assumpF} iii) that $$\phi'(t_0)=w_t(x_0,t_0)\le0.$$
 Thus,  inequality
$(\ref{viscosity})$ is satisfied in the viscosity sense.
\bigskip

\noindent\textbf{Step 4: Conclusion}\\
We deduce that $M(t)\le M(0)$ for all $t\ge 0$, i.e. for all $ (x,t)\in Q$, we have  $$w(x,t)=\frac{|Du(x,t)|^2}{2}\le M(t)\le M(0)=
\max_{{\mathbb{T}}^d}{\frac{|Du(\cdot,0)|^2}{2}};$$
\noindent which ends the proof of Theorem \ref{general case}.
\end{Proofc}

\begin{remark}(\textbf{{Reformulation of (\ref{A}) for quasilinear problems}})

 \item [ i)] If $F$ is in  quasilinear divergence form, independent from $x$, (for an equation $u_t=div(a(Du,u,t))$) i.e. for $F$ as follows
$$F(X,p,u,x,t)= D_pa(p,u,t):X+D_ua(p,u,t)\cdot p$$ for some $C^2$-vector field $a$, then (\ref{A}) is equivalent to
\begin{equation}\label{quasidiv}
-D_pa:X^2+|p|^2D^2_{up}a:X+|p|^2D^2_{uu} a\cdot p\le 0.
\end{equation}
 \item [ ii)] If $F$ is quasilinear but not in divergence form, independent from $x$, i.e.
 $$F(X,p,u,x,t)=A(p,u):X+B(p,u)$$ for some  $C^1$-matrix $A$ and  $C^1$-function $B$,
then (\ref{A}) is equivalent to
\begin{equation}
-A:X^2+|p|^2D_uA:X+|p|^2 D_uB\le 0.
\end{equation}
\end{remark}

\begin{Proofc}{\textbf{Proof of Corollary \ref{PMEGresult}.}}

\noindent Corollary \ref{PMEGresult} follows from the application of Theorem
 \ref{general case}, once we check Assumption \ref{assumpF} and condition  (\ref{constantsolutions}). Condition (\ref{constantsolutions})
is straightforward. We now check Assumption \ref{assumpF}.
\bigskip

\noindent\textbf{Step 1: Checking Assumption \ref{assumpF}  i) and ii)}\\
Equation (\ref{PMEGprob}) is a particular case of (\ref{pborigine}) for  $$F(X,p,u,x,t):=G'(u) \,tr(X)+G''(u)|p|^{2}.$$
As $G\in C^{3}$ and $G'\ge 0$, we get that $F$ satisfies Assumption \ref{assumpF} i) and ii).
\bigskip

\noindent\textbf{Step 2: Checking Assumption \ref{assumpF} iii) }\\
Using (\ref{quasidiv}) with $a(u,p,t)=G'(u)p$ inequality (\ref{A})  can be written as
\begin{equation}\label{assumpG}
-G'(u)\, tr(X^2)+|p^2|G''(u)\,tr(X)+|p|^4G'''(u)\leq 0.
\end{equation}
We have to check   (\ref{assumpG})
for all $(X,p,u,x,t)\in E$ such that $ X\cdot p=0$.
\bigskip

\noindent\textbf{Step 2.1: A preliminary result}\\
We prove that \begin{equation}\label{G^{(3)}}
G'''\le 0\quad \mbox{ on } [\underline{u}, \overline{u}].
\end{equation}
If $G'(u)> 0$, by hypothesis (\ref{Gcondition}), $G^{(3)}(u)\leq 0$.
Now, if $G'(u)= 0$, we consider two cases: \\
\noindent\textbf{\underline{Case 1.}}  Suppose that there exists a sequence $(u_k)_k$ in $[\underline u, \overline u]$ converging to $u$ such that $G'(u_k)>0$,
then $G^{(3)}(u_k)\le 0$, and by continuity of $G^{(3)}$, we get $G^{(3)}(u)\le 0$.\\
\noindent\textbf{\underline{Case 2.} } Suppose that $G'(u)=0$ in   a neighborhood of $u$, then, using the regularity $C^3$ of $G$, $G''(u)=G^{(3)}(u)=0$ in this neighborhood. \\
 In both cases $G^{(3)}(u)\le 0$.
 \bigskip

\noindent\textbf{Step 2.2: The core of the analysis}\\
For $p=0$, inequality (\ref{assumpG}) is satisfied since $G'(u)\ge 0$. We now consider  $p\ne 0:$\\
\noindent\textbf{\underline{Case 1.} }: $d=1$\\
The identity $X\cdot p=0$ implies that $X=0$, hence condition (\ref{assumpG}) is
satisfied provided that $G^{(3)}(u)\leq 0$, which is implied by (\ref{G^{(3)}}).

\noindent\textbf{\underline{Case 2.} }: $d \ge 2$\\
\noindent\textbf{\underline{Subcase 2.1.} }:  $G'(u)=0$\\
By hypothesis (\ref{Gcondition}), we have  $G''(u)=0$, and using (\ref{G^{(3)}}), we have $G^{(3)}(u)\le 0$.
This implies (\ref{assumpG}).

\noindent\textbf{\underline{Subcase 2.2.} }: $G'(u)> 0$\\
\noindent We  set
\begin{equation}\label{abc}
a=G'(u) ,\, -b=|p|^2G''(u)
\mbox { and }-c=|p| ^{4}G^{(3)}(u).\end{equation}
For any $(X,p)\in S^d\times\overline{B(0,M)}$ such that  $X\cdot p=0$, up to change the orthonormal basis of ${\mathbb{R}}^{d}$,
 $X$ can be written as  a diagonal matrix  of eigenvalues $0,{\lambda}_{1},..,{\lambda}_{d-1}$ in the
 direct decomposition $p\oplus p^{\perp}$.
  Writing $$\Lambda=\left( 0,{\lambda}_{1},..., {\lambda}_{d-1}\right)\mbox{ and }
   e=\left(0,1,...,1\right),\,$$
then (\ref{assumpG}) can  be rewritten as

\begin{equation}\label{assumpGeq}
c+be\cdot\Lambda+a{\Lambda}^{2}\ge0.
\end{equation}

\noindent Since
$$|be\cdot\Lambda|=\left|\frac{be}{\sqrt{2a}}\cdot\sqrt{2a}\Lambda\right|
 \leq\frac{1}{2}\left(2a{\Lambda}^{2}+\frac{|be|^{2}}{2a}\right)\leq
  a{\Lambda}^{2}+\frac{|be|^{2}}{4a}, $$
we obtain that (\ref{assumpGeq}) is true provided that $ (d-1)b^{2}\le 4ac$,
 which is given by the hypothesis (\ref{Gcondition}) on $G$.
\end{Proofc}

\begin{Proofc}{\textbf{Proof of Corollary \ref{PMEresult}.}}

\noindent For the existence and uniqueness of the weak solution $u$,   we refer to \cite[Theorem 9.25, p. 218]{vasquez}.
This solution  is constructed as the limit in $L^{\infty}(0,+\infty ; L^{1}({\mathbb{T}}^d))$ of a sequence $u_\epsilon $ in $L^{\infty}(0,+\infty ;W^{1,\infty}({\mathbb{T}}^d))\cap C([0,\infty);L^{1}({\mathbb{T}}^d))$ of smooth positive solutions
 satisfying the same equation with initial data $u_{0\epsilon}=\epsilon+\rho_\epsilon \ast u_0$ for some mollifier $\rho_\epsilon$ with $u_0\ge 0$,
  which implies that

 $$\epsilon+\min_{{\mathbb{T}}^d}u_0\leq u_{0\epsilon}\leq\epsilon+\max_{{\mathbb{T}}^d}u_0\quad\mbox{ and }\quad
 \|Du_{0\epsilon}\|_{L^{\infty}({\mathbb{T}}^d)}\leq\|Du_{0}\|_{L^{\infty}({\mathbb{T}}^d)}.$$ For $1\le m\le1+\frac{4}{3+d}$,
  the function $G(u):=u^{m}$ satisfies conditions of Corollary \ref{PMEGresult} with $\underline{u}=\epsilon+\underset{{\mathbb{T}}^d}{\min }\,u_0$ and $\bar{u}=\epsilon+\underset{{\mathbb{T}}^d}{\max }\, u_0$. Because the solutions are known to be smooth for any $T>0$, there exists $M=M(T)$ such that for all $t\in [0,T]$, we have $\|Du(\cdot,t)\|_{{\mathbb{T}}^d}\leq M.$ Hence, applying a version of Corollary \ref{PMEGresult} for finite time interval $(0,T)$, we obtain for  $t\in [0,T] $
$$\epsilon+\min_{{\mathbb{T}}^d}u_0\leq u_{\epsilon}(\cdot,t)\leq \epsilon+\max_{{\mathbb{T}}^d}u_{0}$$ and $$\|Du_{\epsilon}(\cdot,t)\|_{L^{\infty}({\mathbb{T}}^d)}\leq \|Du_{0\epsilon}\|_{L^{\infty}({\mathbb{T}}^d)}\le\|Du_{0}\|_{L^{\infty}({\mathbb{T}}^d)}.$$
Because $T>0$ is arbitrary,  we recover the bound for all time $t\ge 0$.

\noindent Now $u_\epsilon \rightarrow u$ as $\epsilon$ goes to zero and we recover the expected bounds (\ref{uniformboundonu}) and (\ref{uniformboundonnabla}) for $u$.\end{Proofc}
\begin{remark}
 Notice that bounds (\ref{uniformboundonu}) and (\ref{uniformboundonnabla}) can also be deduced from Corollary \ref{Tphysicalmodele} below, for smooth solutions of the doubly nonlinear diffusion equation
  \begin{equation}\label{laplacienpdimqlq}
  u_t=\Delta_p(u^{m}), \quad (x,t)\in Q.
  \end{equation}
  This works for $p\ge1$,  satisfying $0\le (m-1)(p-1)\le \frac{4}{3+d}$ and where the operator $\Delta_p v$ is defined by $\Delta_p v=\textrm{div} (|\nabla v|^{p-2}\nabla v) $.
 \end{remark}

\section{Applications to models  in hydrology}
\noindent In this section, we apply Theorem \ref{general case} to  hydrological models (\ref{modelehydrologie})
and (\ref{modelehydrologiesimplifie}). To this end,  we will first prove a priori bounds on gradient of  solutions of Problem
(\ref{psimodelegeneral}). We begin
by assumptions on the function $\psi$
\begin{assumption}\label{assumppsi}\textbf{ }\\
Let $d\ge 1$,  $\underline u\le\overline u$, $L>0$ and $\psi:=\psi(u,s)$ be a real function
satisfying the following conditions
\noindent \item [i)]\textbf{\textit{Regularity}}:  $ \psi\in C^{2}([\underline u,\overline u]\times [0,L])$.
 \noindent \item [ii)]\textbf{\textit{Degenerate parabolicity}}:
\noindent $\psi$ satisfies conditions
\begin{eqnarray}
 \psi\ge0 && \mbox{ on } [\underline u,\overline u]\times [0,L],\label{psi+}\\
\psi+2sD_{s}\psi\ge0&&  \mbox{ on }[\underline u,\overline u]\times [0,L].\label{psi<}
\end{eqnarray}
\noindent \item [iii)]\textbf{\textit{Differential inequality}}:
\begin{equation}\label{$B$}
\frac {(d-1)}4(D_{u}\psi)^{2}\le
-\psi D^{2}_{uu}\psi
\quad \mbox{ on }  [\underline u,\overline u]\times [0,L].
\end{equation}
\end{assumption}

\begin{remark}
Hypothesis (\ref{$B$}) on $\psi$ is similar to hypothesis (\ref{Gcondition}) on $G$.
\end{remark}

\begin{corollary}\label{Tphysicalmodele}({\textbf{A  priori gradient bounds for  solutions of Problem
(\ref{psimodelegeneral})}})\\
Suppose that $\psi$ satisfies  Assumption \ref{assumppsi}.
\noindent Let us assume the existence of a function $u\in C^{3}(\overline{Q})$ solution of Problem (\ref{psimodelegeneral}), (\ref{initdata}) satisfying condition (\ref{$H_{u}$}) for $M=\sqrt L$. Then
for all $t\ge0$ we have  a priori bounds (\ref{uniformboundonu}) and (\ref{uniformboundonnabla}).
\end{corollary}

\begin{Proofc}{\textbf{Proof of Corollary \ref{Tphysicalmodele}.}}

\noindent We need only to check Assumption \ref{assumpF} and condition (\ref{constantsolutions}). Condition (\ref{constantsolutions}) is  straightforward.
\bigskip

\noindent\textbf{Step 1: Checking Assumption \ref{assumpF}  i) and ii)}\\
Equation (\ref{PMEGprob}) is a particular case of (\ref{pborigine}) for $$F(X,p,u,x,t):=(\psi(u,|p|^2) Id+
2D_s\psi(u,|p|^2)p\otimes p):X+D_u\psi(u,|p|^2)|p|^{2}.$$
As $\psi \in C^{2}([\underline{u},\overline u]\times[0,L])$ and satisfies (\ref{psi+}) and (\ref{psi<}), we see that
$F$ satisfies Assumption \ref{assumpF}  i) and ii).
\bigskip

\noindent\textbf{Step 2: Checking Assumption \ref{assumpF} iii) }\\
Inequality (\ref{A}) of Assumption \ref{assumpF} iii) can be written as
\begin{equation}\label{assumpFonpsi*}
-(\psi Id+2D_s\psi\, p\otimes p):X^{2}+|p|^2(D_{u}\psi Id+2D^{2}_{us}\psi(u,|p|^2)p\otimes p):X+|p|^4 D ^{2}_{uu}{\psi}\le0.
\end{equation}
for all $(X,p,u,x,t)\in E$ such that $ X\cdot p=0$. This is equivalent to
\begin{equation}\label{assumpFonpsi}
-\psi \,tr(X^{2})+|p|^2D_{u}\psi \,tr(X)+|p|^4 D ^{2}_{uu}{\psi}\le0,
\end{equation}
since $ X\cdot p=0$.
In order to check Assumption \ref{assumpF} iii) we have to check (\ref{assumpFonpsi})
 for all $(X,p,u,x,t)\in E$ such that $ X\cdot p=0$.

\noindent Setting now $$a=\psi,\quad-b=|p|^2D_u\psi\mbox{ and }-c=|p|^2D^{2}_{uu}\psi,$$
inequality (\ref{assumpFonpsi}) reads as \begin{equation}\label{asumponpsiabc} c+b\, tr(X)+a \,tr(X^2)\ge 0,\end{equation}
that we have to check for all $(X,p,u,x,t)\in E$ such that $ X\cdot p=0$.
In  Step 2.1  and Step 2.2 of the proof Corollary \ref{PMEGresult}, we can replace in (\ref{abc}) the function
$G'(u)$ by the function $\psi^s(u):=\psi (u, s)$, for each given  $s=|p|^2$. Then the reasoning there applies here without any change
and shows that  (\ref{asumponpsiabc}) holds true under Assumption \ref{assumpFonpsi} iii).
\end{Proofc}

\begin{corollary}\label{interfaceappliq}(\textbf{Application to model (\ref{modelehydrologie})})\\
\noindent Let \begin{equation}\label{inter} d\ge 1,\quad \delta=({2+2d})^{-\frac12},\quad \underline{u}=\frac 12-\delta,\quad \overline u=\frac 12+\delta\quad\mbox{ and }\quad
 M=1.\end{equation}
Let us assume the existence
 of  $u\in C^{3}(\overline{Q})$, solution  of Problem (\ref{modelehydrologie}), (\ref{initdata}),
such that $u$ and $Du$ satisfy condition (\ref{$H_{u}$}). Then for all $t\ge0$, $u$ and $Du$ satisfy the a priori bounds (\ref{uniformboundonu}) and (\ref{uniformboundonnabla}).
 \end{corollary}

\begin{remark}({\textbf{Weak radial solutions of Problem (\ref{modelehydrologie})}})\\
Assume (\ref{inter}) and let $B(0, R)$ be a ball of radius $R>0.$ Then we have, by \cite{hjm} for $d>1$ and by \cite{DH} for $d=1$, existence and uniqueness of weak radial solutions of Problem  (\ref{modelehydrologie}), (\ref{initdata}) with Neumann boundary conditions. Moreover the solutions satisfy uniform bounds similar to (\ref{uniformboundonu}) and (\ref{uniformboundonnabla}), but on the ball $B(0,R)$.\\
Therefore Corollary \ref{interfaceappliq} appears to be a kind of extension of this result to the non radial case, and for smooth solutions.
 \\Finally let us mention an existence and uniqueness result in \cite{hjm} for radial nondecreasing solutions with $u(x,t)\in [\frac12,1]$.
\end{remark}

\begin{Proofc}{\textbf{Proof of Corollary \ref{interfaceappliq}.}}

\noindent We simply apply Corollary \ref{Tphysicalmodele} with  $\displaystyle{\psi(u, |p|^2)=h(u)f(|p|^2)}$, $h(u)=u(1-u)$ and $f(s)=1/(1+s)$. To this end, we have to check Assumption \ref{assumppsi}.

\noindent Assumption \ref{assumppsi} i) and ii) are satisfied for
$0\le\underline{u}\le \bar{u} \le1$ and $L=1=M^2$. Indeed $$\psi(s)+2sD_s\psi(s)=h(u)(f(s)+2sf'(s))=h(u)\frac{(1-s)}{(1+s)^2}\ge0,$$
for all $(u,s)\in[\underline u,\overline u]\times[0,L]$.

\noindent Inequality (\ref{psi<}) of Assumption \ref{assumppsi} iii) can be written as
\begin{equation}\label{condition}
\frac{(d-1)}{4}\left(\frac{(h'(u))^{2}}{-h(u)\,h''(u)}\right)
\le 1.\end{equation}
Setting $u=\frac12+v$, (\ref{condition}) is equivalent to
\begin{equation}\label{func(u)}\frac{(d-1)}{2}\leq \frac{\frac14-v^2}{v^2},\end{equation}
which means that $|v|\le {(2+2d)^{-\frac12}}=\delta$, i.e. $u\in[\underline u,\overline u]$.
This shows that Assumption \ref{assumppsi} iii) holds true.
\end{Proofc}

\noindent We also have
\begin{corollary}(\textbf{{Application to model (\ref{modelehydrologiesimplifie})}})\label{coromodelinterfacesimpl}\\
\noindent Let $d\ge 1$, $\delta=({2+2d})^{-\frac12}$, $\underline{u}=\frac 12-\delta$, $\overline u=\frac 12+\delta$,
and $M>0$ be a real number. Let us assume the existence
 of  $u\in C^{3}(\overline{Q})$, solution  of Problem (\ref{modelehydrologiesimplifie}), (\ref{initdata}),
such that $u$ and $Du$ satisfy condition (\ref{$H_{u}$}).
 Then for all $t\ge0$, $u$ and $Du$ satisfy  a priori bounds (\ref{uniformboundonu}) and (\ref{uniformboundonnabla}).
 \end{corollary}
The proof of  Corollary \ref{coromodelinterfacesimpl} is similar  to the one  of Corollary \ref{interfaceappliq},
taking the same function $h(u)=u(1-u)$ but with $f(s)=1$. Corollary \ref{coromodelinterfacesimpl} can also be obtained
from Corollary
\ref{PMEGresult}.


\section{Appendix}

\begin{Proofc}{\textbf{Proof of Proposition \ref{comparisononu}.}}

\noindent Let $z:=u-v$,
then we have $$\begin{array}{ll}
z_t=&F(D^2u,Du,u,x,t)-F(D^2v,Du,u,x,t)\\
&\\
&+F(D^2v,Du,u,x,t)-F(D^2v,Dv,u,x,t)\\
&\\
&+F(D^2v,Dv,u,x,t)-F(D^2v,Dv,v,x,t).
\end{array}
$$
By mean of Taylor expansion, $z$ satisfies the following problem
 $$\left\{\begin{array}{ll}
  z_{t}=A(x,t):D^2 z+B(x,t)\cdot Dz+C(x,t)z,& (x,t) \in {\mathbb{T}}^d\times (0,T),\\
  &
               \\
     z(x,0)=u(x,0)- v(x,0), & x \in {\mathbb{T}}^d,
   \end{array}\right.$$
   where
$$\left\{\begin{array}{l}
A(x,t)=\int_{0}^{1}D_{X}F(\theta D^2u+(1-\theta)D^2v,Du,u,x,t) d\theta, \\
\\
B(x,t)=\int_{0}^{1}D_{p}F(D^2v,\theta Du+(1-\theta)Dv,u,x,t) d\theta,\\
\\
C(x,t)=\int_{0}^{1}D_{u}F(D^2v,Dv,\theta u+(1-\theta)v,x,t) d\theta.
\end{array}
\right.$$
Note that Assumption \ref{assumpF} ii) implies that $D_{X}F\ge 0$ and then $A\ge 0$.
\noindent Now, since $z(\cdot,0)\leq0$, the standard maximum
principle implies that $z(\cdot,t)\leq0$ for all $t\in[0,T]$
(see, for instance, \cite[Theorem 9, p. 369]{Evans} applied to $w:=e^{-\lambda t}z$ for a large $\lambda>0$).

\end{Proofc}
\begin{center}\textsc{Acknowledgements}\end{center}
This work has been financially supported  by the Lebanese Association for Scientific Research and the Project AUF-MERSI.

\end{document}